\input amstex\documentstyle{amsppt}  
\pagewidth{12.5cm}\pageheight{19cm}\magnification\magstep1
\NoBlackBoxes
\topmatter
\title On the Steinberg character of a semisimple $p$-adic group\endtitle
\author Ju-Lee Kim and George Lusztig\endauthor
\address{Department of Mathematics, M.I.T., Cambridge, MA 02139}\endaddress
\dedicatory{Dedicated to Robert Steinberg on the occasion of his 90-th birthday}\enddedicatory
\thanks{Both authors are supported in part by the National Science Foundation}\endthanks
\endtopmatter   
\document

\define\dw{\dot w}

\define\uw{\un w}

\define\uK{\un K}

\define\si{\sim}

\define\sqc{\sqcup}

\redefine\spa{\spadesuit}
\define\part{\partial}
\define\em{\emptyset}

\define\ra{\rangle}
\define\n{\notin}
\define\iy{\infty}
\define\m{\mapsto}
\define\do{\dots}
\define\la{\langle}
\define\bsl{\backslash}

\define\sub{\subset}    

\define\T{\times}

\define\nl{\newline}
\redefine\i{^{-1}}

\define\un{\underline}

\define\ot{\otimes}

\define\Ad{\text{\rm Ad}}

\define\Aut{\text{\rm Aut}}

\define\tr{\text{\rm tr}}

\define\a{\alpha}
\redefine\b{\beta}
\redefine\c{\chi}
\define\g{\gamma}
\redefine\d{\delta}

\define\ph{\phi}

\define\r{\rho}
\define\s{\sigma}
\redefine\t{\tau}

\redefine\l{\lambda}

\define\Om{\Omega}

\define\CC{\bold C}

\define\NN{\bold N}

\define\RR{\bold R}
\define\SS{\bold S}

\define\ZZ{\bold Z}

\define\cb{\Cal B}
\define\cc{\Cal C}

\define\cl{\Cal L}

\define\co{\Cal O}

\define\cw{\Cal W}

\define\cx{\Cal X}
\define\cy{\Cal Y}

\define\fg{\frak g}

\define\fl{\frak l}

\define\fn{\frak n}

\define\fp{\frak p}

\define\ft{\frak t}

\define\fH{\frak H}

\define\fT{\frak T}

\define\sha{\sharp}

\head 1. Introduction\endhead
\subhead 1.1\endsubhead
Let $K$ be a nonarchimedean local field and let $\uK$ be a maximal unramified field extension of $K$. 
Let $\co$ (resp. $\un\co$) be the ring of integers of $K$ (resp. $\uK$) and let $\fp$ (resp. $\un\fp$)
be the maximal ideal of $\co$ (resp. $\un\co$). Let $\uK^*=\uK-\{0\}$.
We write $\co/\fp=F_q$, a finite field with $q$ elements, of characteristic $p$. 

Let $G$ be a semisimple almost simple algebraic group defined and split over $K$ with a given 
$\co$-structure compatible with the $K$-structure. 

If $V$ is an admissible representation of $G(K)$ of finite length, we denote by $\ph_V$ the character of $V$
in the sense of Harish-Chandra, viewed as a $\CC$-valued function on the set $G(K)_{rs}:=G_{rs}\cap G(K)$.
(Here $G_{rs}$ is the set of regular semisimple elements of $G$ and $\CC$ is the field of complex numbers.) 

In this paper we study the restriction of the function $\ph_V$ to:

(a) a certain subset $G(K)_{vr}$ of $G(K)_{rs}$, that is to the set of very regular elements in $G(K)$ (see 
1.2), in the case where $V$ is the Steinberg representation of $G(K)$ and

(b) a certain subset $G(K)_{svr}$ of $G(K)_{vr}$, that is to the set of split very regular elements in 
$G(K)$ (see 1.2), in the case where $V$ is an 
irreducible admissible representation of $G(K)$ with nonzero vectors fixed by an Iwahori subgroup.

In case (a) we show that $\ph_V(g)$ with $g\in G(K)_{rs}$ is of the form $\pm q^n$ with 
$n\in\{0,-1,-2,\do\}$ (see Corollary 3.4) with more precise information when $g\in G(K)_{svr}$ (see 
Theorem 2.2) or when $g\in G(K)_{cvr}$ (see Theorem 3.2); in case (b) we show (with some restriction on 
characteristic) that $\ph_V(g)$ with $G(K)_{svr}$ can be expressed as a trace of a certain element of an 
affine Hecke algebra in an irreducible module (see Theorem 4.3).

Note that the Steinberg representation $\SS$ is an irreducible admissible representation of $G(K)$ with 
a one dimensional subspace invariant under an Iwahori subgroup on which the affine Hecke algebra acts 
through the ``sign'' representation, see \cite{MA}, \cite{S}. This is a $p$-adic analogue of the 
Steinberg representation \cite{St} of a reductive group over $F_q$. 
In \cite{R}, it is proved that $\ph_\SS(g)\ne0$ for any $g\in G(K)_{rs}$.

\subhead 1.2\endsubhead
Let $g\in G_{rs}\cap G(\uK)$. Let $T'=T'_g$ be the maximal torus of $G$
that contains $g$. We say that $g$ is very regular (resp. compact very regular)
if $T'$ is split over $\uK$ and for any root $\a$ with respect to $T'$
viewed as a homomorphism $T'(\uK)@>>>\uK^*$ we have 

$\a(g)\n(1+\un\fp)$ (resp. $\a(g)\in\un\co$, $\a(g)\n(1+\un\fp)$).

Let $G(\uK)_{vr}$ (resp. $G(\uK)_{cvr}$) be the set of elements in $G(\uK)$ which are very regular
(resp. compact very regular). We write $G(K)_{vr}=G(\uK)_{vr}\cap G(K)$, $G(K)_{cvr}=G(\uK)_{cvr}\cap G(K)$.
Let $G(K)_{svr}$ be the set of all $g\in G(K)_{vr}$ such that $T'_g$ is split over $K$.

\subhead 1.3\endsubhead {\it Notation.} 
Let $K^*=K-\{0\}$ and let $v:K^*@>>>\ZZ$ be the unique (surjective) homomorphism
such that $v(\fp^n-\fp^{n+1})=n$ for any $n\in\NN$. For $a\in K^*$ we set $|a|=q^{-v(a)}$. 

We fix a maximal torus $T$ of $G$ defined and split over $K$. 
Let $Y$ (resp. $X$) be the group of cocharacters (resp. characters) of the algebraic group $T$. 
Let $\la,\ra:Y\T X@>>>\ZZ$ be the obvious pairing.
Let $R\sub X$ be the set of roots of $G$ with respect to $T$, let $R^+$ be a set of positive roots for $R$
and let $\Pi$ be the set of simple roots of $R$ determined by $R^+$. We write $\Pi=\{\a_i;i\in I_0\}$. 
Let $R^-=R-R^+$.
Let $Y^+$ (resp. $Y^{++}$)  be the set of all $y\in Y$ such that 
$\la y,\a\ra\ge0$ (resp. $\la y,\a\ra>0$) for all $\a\in R^+$
We define $2\r\in X$ by $2\r=\sum_{a\in R^+}\a$.

We have canonically $T(K)=K^*\ot Y$; we define a homomorphism $\c:T(K)@>>>Y$ 
by $\c(\l\ot y)=v(\l)y$ for any $\l\in K^*,y\in Y$. For any $y\in Y$ we set $T(K)_y=\c\i(y)$.
For $y\in Y$ let $T(K)^\spa_y=T(K)_y\cap G(K)_{svr}$. 
Note that if $y\in Y^{++}$ then $T(K)^\spa_y=T(K)_y$.

For each $\a\in R$ let $U_\a$ be the corresponding root subgroup of $G$.

Let $G(K)'$ be the derived subgroup of $G(K)$. 

\head 2. Calculation of $\ph_\SS$ on $G_(K)_{svr}$\endhead
\subhead 2.1\endsubhead
Let $\cw\sub\Aut(T)$ be the Weyl group of $G$ regarded as a Coxeter group;
for $i\in I_0$ let $s_i$ be the simple reflection in $\cw$ determined by $\a_i$. We can also view $\cw$ as 
a subgroup of $\Aut(Y)$ or $\Aut(X)$. Let $w=w_0$ be the longest element of $\cw$. 
For any $J\sub I_0$ let $\cw_J$ be the subgroup of $\cw$ generated by 
$\{s_i;i\in J\}$ and let $R_J$ be the set of $\a\in R$ such that $\a=w(\a_i)$ for some $w\in\cw_J,i\in J$. 
Let $R_J^+=R_J\cap R^+$, $R_J^-=R_J-R_J^+$. 

 Let $\fg$ be the Lie algebra of $G$; let $\ft\sub\fg$ be the Lie algebra of $T$. 
For any $J\sub I_0$ let $\fl_J$ be the Lie subalgebra of $\fg$ spanned by $\ft$ and by the root spaces 
corresponding to roots in $R_J$; let $\fn_J$ be 
the Lie subalgebra of $\fg$ spanned by the root spaces corresponding to roots in $R^+-R_J^+$.

According to \cite{C1}, $\ph$ is an alternating sum of characters of representations induced from one 
dimensional representations of various parabolic subgroups of $G$ defined over $K$. 
From this one can deduce that, if 
$t\in T(K)\cap G(K)_{rs}$, then
$$\ph_\SS(t)=\sum_{J\sub I}(-1)^{\sha J}\sum_{w\in{}^J\cw}\d_J(w(t))^{1/2}D_{I,J}(w(t))^{-1/2}$$
where for any $J\sub I$ and $t'\in T(K)\cap G(K)_{rs}$ we set
$$D_{I,J}(t')=|\det(1-\Ad(t')|_{\fg/\fl_J})|,$$
$$\d_J(t')=|\det(\Ad(t')|_{\fn_J})|,$$
and ${}^J\cw$ is a set of representatives for the cosets $\cw_J\bsl\cw$. (It will be convenient to assume 
that ${}^J\cw$ is the set of representatives of minimal length for the cosets $\cw_J\bsl\cw$.) 
Here for a real number $a\ge0$ we denote by $a^{1/2}$ or $\sqrt{a}$ the $\ge0$ square root of $a$.
We have the following result. (We write $\ph$ instead of $\ph_\SS$.)

\proclaim{Theorem 2.2} Let $y\in Y^+$ and let $t\in T(K)^\spa_y$. Then $\ph(t)=q^{-\la y,2\r\ra}$.
\endproclaim

\subhead 2.3\endsubhead
More generally let $t\in T(K)^\spa_y$ where $y\in Y$. By a standard property of Weyl chambers 
there exists $w\in\cw$ such that $w(y)\in Y^+$. Let $t_1=w(t)$. Then the theorem is applicable to
$t_1$ and we have $\ph(t)=\ph(t_1)=q^{-\la w(y),2\r\ra}$.

\subhead 2.4\endsubhead
Let $y'=w_0(y), t'=w_0(t)$. We have $\ph_\SS(t)=\ph_\SS(t')$, $t'\in T(K)^\spa_{y'},-y'\in Y^+$. We show:

(a) if $\b\in R^+$ then $v(1-\b(t')))=v((\b(t'))$; if $\b\in R^-$ then $v(1-\b(t')))=0$.
\nl
Assume first that $\b\in R^+$. 
If $v(\b(t'))\ne0$ then $v(\b(t'))<0$ (since $\la y',\b\ra\ne0$, $\la y',\b\ra\le0$)
hence $v(1-\b(t')))=v((\b(t'))$. If $v(\b(t'))=0$ then $\b(t')-1\in\co-\fp$ hence
$v(1-\b(t')))=0=v((\b(t'))$ as required. 

Assume next that $\b\in R^-$.
If $v(\b(t'))\ne0$ then $v(\b(t'))>0$ (since $\la y',\b\ra\ne0$, $\la y',\b\ra\ge0$)
hence $v(1-\b(t')))=0$. If $v(\b(t'))=0$ then $\b(t')-1\in\co-\fp$ hence $v(1-\b(t')))=0$ as required. 
 
For any $w\in\cw,J\sub I$ we have:
$$\align&D_{I,J}(w(t'))=\prod_{\a\in R-R_J}q^{-v(1-\a(w(t')))}\\&=\prod_{\a\in R-R_J;w\i\a\in R^+}
q^{-v(\a(w(t')))}=\prod_{\a\in R-R_J;w\i\a\in R^+}q^{-\la y',w\i\a\ra},\endalign$$
$$\d_J(w(t'))=\prod_{\a\in R^+-R^+_J}q^{-v(\a(w(t')))}=\prod_{\a\in R^+-R^+_J}q^{-\la y',w\i\a\ra},$$
$$D_I(t')=\prod_{\a\in R^+}q^{-\la y',\a\ra}.$$
(We have used (a) with $\b=w\i(\a)$.) We see that
$$\ph(t)=\ph(t')=\sum_{J\sub I}(-1)^{\sha J}\sum_{w\in{}^J\cw}\sqrt{q}^{-\la y',x_{w,J}\ra}$$
where for $w\in{}^J\cw$ we have
$$\align&x_{w,J}
=\sum_{\a\in R^+-R^+_J}w\i\a-\sum_{\a\in R-R_J;w\i\a\in R^+}w\i\a\\&
=\sum_{\a\in R^+-R^+_J;w\i(\a)\in R^-}w\i\a-\sum_{\a\in R^--R^-_J;w\i(\a)\in R^+}w\i\a\\&
=2\sum_{\a\in R^+-R^+_J;w\i\a\in R^-}w\i\a\in X.\endalign$$
For $w\in{}^J\cw$ we have $\a\in R^+_J\implies w\i\a\in R^+$ hence 
$$\sum_{\a\in R^+-R^+_J;w\i\a\in R^-}w\i\a=\sum_{\a\in R^+;w\i\a\in R^-}w\i\a$$
so that $x_{w,J}=x_w$ where 
$$x_w=2\sum_{\a\in R^+;w\i\a\in R^-}w\i\a\in X.$$
Thus we have
$$\ph(t)=\sum_{J\sub I}(-1)^{\sha J}\sum_{w\in{}^J\cw}\sqrt{q}^{-\la y', x_w\ra}
=\sum_{w\in\cw}c_w\sqrt{q}^{-\la y',x_w\ra}$$
where for $w\in\cw$ we set
$$c_w=\sum_{J\sub I;w\in{}^J\cw}(-1)^{\sha J}.$$
For $w\in\cw$ let $\cl(w)=\{i\in I;s_iw>w\}$ where $<$ is the standard partial order on $\cw$. 
For $J\sub I$ we have $w\in{}^J\cw$ if and only if $J\sub\cl(w)$. Thus,
$$c_w=\sum_{J\sub\cl(w)}(-1)^{\sha J}$$
and this is $0$ unless $\cl(w)=\em$ (that is $w=w_0$) when $c_w=1$. Note also that
$x_{w_0}=-4\r$. Thus we have 
$$\ph(t)=c_{w_0}\sqrt{q}^{-\la y',x_{w_0}\ra}=q^{\la y',2\r\ra}=q^{-\la y,2\r\ra}.$$
Theorem 2.2 is proved.

\subhead 2.5\endsubhead
Assume now that $\t\in T(K)$ satisfies the following condition: for any $\a\in R$ we have 
$\a(\t)-1\in\fp-\{0\}$ so that $\a(\t)-1\in\fp^{n_\a}-\fp^{n_\a+1}$ for a well defined integer 
$n_\a\ge1$. Note that $n_{-\a}=n_\a$ and $v(1-\a(\t))=n_\a\ge1$ for all $\a\in R$. Hence 
$$\ph(\t)=\sum_{J\sub I}(-1)^{\sha J}\sum_{w\in{}^J\cw}
q^{\sum_{\a\in R}n_\a/2-\sum_{\a\in R_J}n_{w\i(\a)}/2}.$$
Thus,
$$\ph(\t)=\sha(\cw)q^{\sum_{\a\in R}n_\a/2}+\text{strictly smaller powers of }q.$$
In the case where $K$ is the field of power series over $F_q$, the leading term 

$\sha(W)q^{\sum_{\a\in R}n_\a/2}$
\nl
is equal to $\sha(\cw)q^m$ where
$m$ is the dimension of the ``variety'' of Iwahori subgroups of $G(\uK)$ 
that contain the topologically unipotent element $\t$ (see \cite{KL2}).

\head 3. Calculation of $\ph_\SS$ on $G_(K)_{vr}$\endhead
\subhead 3.1\endsubhead
We will again write $\ph$ instead of $\ph_\SS$. In this section we assume that we are given 
$\g\in G(K)_{vr}$.
Let $T'=T'_\g$. Note that $T'$ is defined over $K$; let $A'$ be the largest $K$-split torus of $T'$.
For any parabolic subgroup $P$ of $G$ defined over $K$ such that $\g\in P$ we set
$\d_P(\g)=|\det(\Ad(\g)|_{\fn})|$ where $\fn$ is the
Lie algebra of the unipotent radical of $P$.

Let $\cx$ be the set of all pairs $(P,A)$ where $P$ is a parabolic subgroup of $G$ defined over $K$ and $A$
is the unique maximal $K$-split torus in the centre of some Levi subgroup of $P$ defined over $K$; then that
Levi subgroup is uniquely determined by $A$ and is denoted by $M_A$.
Let $\cx'=\{(P,A)\in\cx;A\sub A'\}$. According to Harish-Chandra \cite{H} we have
$$\ph(\g)=(-1)^{\dim T}\sum_{(P,A)\in\cx'}(-1)^{\dim A}\d_P(\g)^{1/2}D_{G/M_A}(\g)^{-1/2}\tag a$$
where $D_{G/M_A}(\g)=|\det(1-\Ad(\g)|_{\fg/\fl})|$ (we denote by $\fl$ the Lie algebra of $M_A$).

\proclaim{Theorem 3.2} Assume in addition that $\g\in G(K)_{cvr}$. Then 

$\ph(\g)=(-1)^{\dim T-\dim A'}$.
\endproclaim
From our assumptions we see that for any $(P,A)\in\cx'$ we have $\d_P(\g)=1=D_{G/M_A}(\g)$. Hence 3.1(a) 
becomes
$$\ph(\g)=(-1)^{\dim T}\sum_{(P,A)\in\cx'}(-1)^{\dim A}.$$
Let $\cy$ be the group of cocharacters of $A'$ and let $\fH=\cy\ot\RR$. The real vector space $\fH$ can be
partitioned into facets $F_{P,A}$ indexed by $(P,A)\in\cx'$ such that $F_{P,A}$ is homeomorphic to
$\RR^{\dim A}$. Note that the Euler characteristic with compact support of $F_{P,A}$ is $(-1)^{\dim A}$ and 
the Euler characteristic with compact support of $\fH$ is $(-1)^{\dim_\RR\fH}=(-1)^{\dim A'}$. Using the 
additivity of the Euler characteristic with compact support we see that
$\sum_{(P,A)\in\cx'}(-1)^{\dim A}=(-1)^{\dim A'}$. Thus, $\ph(\g)=(-1)^{\dim T-\dim A'}$, as required. \qed

\subhead 3.3\endsubhead
In the setup of 3.1 let $P_\g$ be the parabolic subgroup of $G$ associated to $\g$ as in \cite{C2}. Note 
that $P_\g$ is defined over $K$. The following result can be deduced by combining Theorem 3.2 with the 
results in \cite{C2} and with Proposition 2 of \cite{R}.

\proclaim{Corollary 3.4} We have $\ph(\g)=(-1)^{\dim T-\dim A'}\d_{P_{\g}}(\g)$.
\endproclaim

\head 4. Iwahori spherical representations: split elements\endhead
\subhead 4.1\endsubhead
Let $B$ be the subgroup of $G(K)$ generated by $U_\a(\co),(\a\in R^+)$, $U_\a(\fp),(\a\in R^-)$ and 
$T(K)_0$. (The subgroups $U_\a(\co),U_\a(\fp)$ of $U_\a$ are defined by the $\co$-structure of $G$.
We have $B\in\cb$ where $\cb$ is the set of Iwahori subgroups of $G(K)$. Note that $B\sub G(K)'$.
For any $\a\in R$ we choose an isomorphism $x_\a:K@>\si>>U_\a(K)$ (the restriction of an isomorphism of 
algebraic groups from the additive group to $U_\a$) which carries $\co$ onto $U_\a(\co)$ and
$\fp$ onto $U_\a(\fp)$.
We set $W:=Y\cdot\cw$ with $Y$ normal in $W$ (recall that $\cw$ acts naturally on $Y$).
Let $Y'$ be the subgroup of $Y$ generated by the coroots. Then $W':=Y'\cdot\cw$ is naturally a subgroup
of $W$. According to \cite{IM}, $W$ is an extended
Coxeter group (the semidirect product of the Coxeter group $W'$ with the finite abelian group $Y/Y'$) with 
length function 
$$l(yw)=\sum_{\a\in R^+;w\i(\a)\in R^+}||\la y,\a\ra||+\sum_{\a\in R^+;w\i(\a)\in R^-}||\la y,\a\ra-1||$$
where $||a||=a$ if $a\ge0$, $||a||=-a$ if $a<0$.
According to \cite{IM}, the set of double cosets $B\bsl G(K)/B$ is in bijection with $W$; to $yw$ (where 
$y\in Y,w\in\cw$) corresponds the double coset $\Om_{yw}$ containing $T(K)_y\dw$ (here $\dw$ is an element in
$G(\co)$ which normalizes $T(K)_0$ and acts on it in the same way as $w$); moreover, 
$\sha(\Om_{yw}/B)=\sha(B\bsl\Om_{yw})=q^{l(yw)}$ for any $y\in Y$, $w\in\cw$. For example, if $y\in Y^{++}$ 
then $l(y)=\la y,2\r\ra$.

Let $H$ be the algebra of $B$-biinvariant functions $G(K)@>>>\CC$ with compact support with respect to 
convolution (we use the Haar measure $dg$ on $G(K)$ for which $vol(B)=1$).
For $y,w$ as above let $\fT_{yw}\in H$ be the characteristic function of $\Om_{yw}$. 
Then the functions $\fT_{\uw}$, $\uw\in W$, form a $\CC$-basis of $H$ and according to \cite{IM} we have

$\fT_{\uw}\fT_{\uw'}=\fT_{\uw\uw'}$ if $\uw,\uw'\in W$ satisfy $l(\uw\uw')=l(\uw)+l(\uw')$,

$(\fT_{\uw}+1)(\fT_{\uw}-q)=0$ if $\uw\in W',l(\uw)=1$.
\nl
In other words, $H$ is what now one calls the Iwahori-Hecke algebra of the (extended) Coxeter group $W$ with 
parameter $q$. 

\subhead 4.2\endsubhead
Let $\cc_0^\iy(G(K))$ be the vector space of locally constant functions with compact support from $G(K)$ to 
$\CC$. 
Let $(V,\s)$ be an irreducible admissible representation of $G(K)$ such that the space $V^B$ of 
$B$-invariant vectors in $V$ is nonzero. If $f\in\cc_0^\iy(G(K))$ then there is a well defined linear map 
$\s_f:V@>>>V$ such that for any $x\in V$ we have $\s_f(x)=\int_Gf(g)\s(g)(x)dg$. This linear map has finite 
rank hence it has a well defined trace $\tr(\s_f)\in\CC$. From the definitions we see that for 
$f,f'\in\cc^\iy_0(G(K))$ we have $\s_{f*f'}=\s_f\s_{f'}:V@>>>V$ where $*$ denotes convolution. If $f\in H$ ,
then $\s_f$ maps $V$ into $V^B$ and $\tr(\s_f)=\tr(\s_f|_{V^B})$. (Recall that $\dim V^B<\iy$.) We see that 
the maps $\s_f|_{V^B}$ define a (unital) $H$-module structure on $V^B$. It is known \cite{BO} that the 
$H$-module $V^B$ is irreducible. Moreover for $\uw\in W$ we have
$\tr(\s_{\fT_{\uw}})=\tr(\fT_{\uw})$
where the trace in the right side is taken in the $H$-module $V^B$. We have the following result.

\proclaim{Theorem 4.3} Assume that $K$ has characteristic zero and that $p$ is
sufficiently large. Let $y\in Y^+$ and let $t\in T(K)^\spa_y$. We have
$$\ph_V(t)=q^{-\la y,2\r\ra}\tr(\fT_y)$$
where the trace in the right side is taken in the irreducible $H$-module $V^B$.
\endproclaim
An equivalent statement is that 
$$\ph_V(t)=\tr(\s_{\fT_y})/vol(\Om_y).$$
(Recall that $\fT_y$ in the right hand side is the characteristic function of $\Om_y=BT(K)_yB$.)

The assumption on characteristic in the theorem is needed only to be able to use a result from \cite{AK}, 
see 5.1$(\dag)$. We expect that the theorem holds without that assumption.

In the case where $y=0$ the theorem becomes:

(a) {\it If $t\in T(K)\cap G_{cvr}$ then $\ph_V(t)=\dim(V^B)$.}
\nl
As pointed out to us by R. Bezrukavnikov and S. Varma, in the special case where $y\in Y^{++}$, Theorem
4.3 can be deduced from results in \cite{C2}.

\subhead 4.4\endsubhead
In the case where $V=\SS$, see 1.1, for any $y\in Y^+$, $\fT_y$ acts on
the one dimensional vector space $V^B$ as the identity map so that $\ph_V(t)=q^{-\la y,2\r\ra}$; 
we thus recover Theorem 2.2 (which holds without assumption on the characteristic).

\head 5. Proof of Theorem 4.3\endhead
\subhead 5.1\endsubhead
Let $B=B_0,B_1,B_2,\cdots$ be the strictly decreasing Moy-Prasad filtration of $B$. In \cite{MP}, this is a 
sequence associated to a point $x$ in the building such that $B=G_{x,0}$. Note that each $B_i/B_{i+1}$ is 
abelian. Let $T_n:=T(K)\cap B_n$. Applying Corollary 12.11 in \cite{AK} to $\ph_V$, we have 

\smallskip

\noindent
$(\dag)$ \qquad\qquad
$\ph_V$ is constant on the $\Ad(G)$-orbit ${}^G\!(tT_1)$ of $tT_1$.

\proclaim{Lemma 5.2} Let $n\ge1$. For any $t'\in T(K)^\spa_y$ and $z\in B_n$, there exist $g\in B_n$, 
$t''\in T_n$ and $z'\in B_{n+1}$ such that $\Ad(g)(t'z)=t't''z'$.
\endproclaim

\noindent
{\it Proof.}
Let $Z=\{\a\in R\mid U_\a\cap B_n\supsetneq U_\a\cap B_{n+1}\}$. If $Z=\em$,
 $B_n=T_nB_{n+1}$. Hence, $z=t''z'$ for some $t''\in T_n$ and $z'\in B_{n+1}$ and one can take $g=1$. 
If $Z\neq\em$, there are $a_\a\in K$, $\a\in Z$ such that $x_\a(a_\a)\in B_n$ and 
$z\equiv\prod_{\a\in Z}x_\a(a_\a)\pmod{T_nB_{n+1}}$. Such $a_\a$ can be chosen independent of the
order of $\prod$ since $B_n/T_nB_{n+1}$ is abelian. Take 
$g=\prod_{\a\in Z}x_\a((1-\a(t'{}^{-1}))^{-1}a_\a)$. 
Then, we have $t'{}^{-1}gt'g^{-1}\equiv z^{-1}\pmod{T_nB_{n+1}}$. 
Moreover, since $y\in Y^+$, we have
$|1-\a(t'{}^{-1})|\ge1$ and thus $g\in B_n$. 
(We argue as in 2.4(a). Assume first that $\a\in R^+$. 
If $v(\a(t'{}\i))\ne0$ then $v(\a(t'{}\i))<0$ (since $\la y,\a\ra\ne0$, $\la y,\a\ra\ge0$)
hence $v(1-\a(t'{}\i)))=v((\a(t'{}\i))<0$ and $|1-\a(t'{}^{-1})|>1$.
If $v(\a(t'{}\i))=0$ then $\a(t'{}\i)-1\in\co-\fp$ 
hence $v(1-\a(t'{}\i)))=0$ and $|1-\a(t'{}^{-1})|=1$ as required. 
Assume next that $\a\in R^-$. 
If $v(\a(t'{}\i))\ne0$ then $v(\a(t'{}\i))>0$ (since $\la y,\a\ra\ne0$, $\la y,\a\ra\le0$)
hence $v(1-\a(t'{}\i)))=0$ and $|1-\a(t'{}^{-1})|=1$ as required. 
If $v(\a(t'{}\i))=0$ then $\a(t'{}\i)-1\in\co-\fp$ 
hence $v(1-\a(t'{}\i)))=0$ and $|1-\a(t'{}^{-1})|=1$ as required.)

Writing $\Ad(g)(t'z)=t'\cdot(t'{}^{-1}gt'g^{-1})\cdot(g zg^{-1})$, we observe that 
$gzg^{-1}\equiv z\pmod{B_{n+1}}$ and $t'{}^{-1}gt'g^{-1}z\in T_nB_{n+1}$. 
Hence $\Ad(g)(t'z)$ can be written as $t't''z'$ with $t''\in T_n$ and $z'\in B_{n+1}$. \hfill\qed 

\proclaim{Lemma 5.3} $B_1tB_1\subset{} ^G\!(tT_1)$.
\endproclaim
\noindent{\it Proof.}
It is enough to show that $tB_1\subset{}^G\!(tT_1)$. Let $t_0z_1\in tB_1$ with $t_0=t$ and 
$z_1\in B_1$. We will construct inductively sequences $g_1,g_2\cdots$, $t_1,t_2\cdots$ and 
$z_1,z_2, \cdots$ such that 
$\Ad(g_k\cdots g_2g_1)(t_0z_1)=\Ad(g_k)(t_0t_1\cdots t_{k-1}z_k)=(t_0t_1\cdots t_k)z_{k+1}$ with 
$g_i\in B_i$, $t_i\in T_i$ and $z_i\in B_i$.

Applying Lemma 5.2 to $n=1$, $t'=t_0$ and $z=z_1$, we find $t_1\in T_1$ and $z_2\in B_2$ such that 
$g_1t_0z_1g_1^{-1}=t_0t_1z_2$ with $t_1\in T_1$ and $z_2\in B_2$. 
Suppose we found $g_i\in B_i$, $z_{i+1}\in B_{i+1}$ and $t_i\in T_i$ for $i=1,\cdots k$ where
$k\ge1$.
Applying Lemma 5.2 to $n=k+1$, $t'=t_0t_1\cdots t_k$ and $z=z_{k+1}$, we find $g_{k+1}\in B_{k+1}$, 
$t_{k+1}\in T_{k+1}$ and $z_{k+2}\in B_{k+2}$ so that 
$g_{k+1}t_0t_1\cdots t_k z_{k+1}g_{k+1}^{-1}=\Ad(g_{k+1}\cdots g_2g_1)(t_0z_1)
=t_0t_1t_2\cdots t_{k+1}z_{k+2}$. (To apply Lemma 5.2 we note that 
$t'\in T(K)^\spa_y$ since $t_0\in T(K)^\spa_y$ and $t_1\cdots t_k\in T_1$ so that for any 
$\a\in R$ we have $\a(t_1\cdots t_k)\in 1+\fp$.)
Taking $g\in B_1$ be the limit of $g_k\cdots g_2g_1$ as $k\rightarrow\infty$, we have 
$\Ad(g)(t_0z_1)\in tT_1$. \hfill\qed

\subhead 5.4\endsubhead
Continuing with the proof of Theorem 4.3, we note that by Lemma 5.3 and 5.1$(\dag)$, for the 
characteristic function $f_t$ of $B_1tB_1$ we have

\noindent$(\ast)$
$$\tr(\s_{f_t})=\int_G f_t(g)\ph_V(g)\,dg=\int_{B_1tB_1}\ph_V(t)\,dg=vol(B_1tB_1)\ph_V(t).$$
Thus it remains to show that
$$\tr(\s_{f_t})/vol(B_1tB_1)=\tr(\s_{\fT_y})/vol(BtB).$$
Since $B_1$ is normalized by $B$, $B$ acts on $V^{B_1}$. Moreover, since $V$ is irreducible and 
$V^B\neq 0$, $B$ acts trivially on $V^{B_1}$ (otherwise, there would exist a nonzero subspace of $V$ on 
which $B$ acts through a nontrivial character of $B/B_1$; since $V^B\neq0$ we see that $(V,\sigma)$ would 
have two distinct cuspidal supports, a contradiction). Thus we have $V^{B_1}=V^B$.
Since $\s_{f_t}$ and $\s_{\fT_y}$ have image contained in $V^{B_1}=V^B$,  it is enough to show that
$$\tr(\s_{f_t}|_{V^B})/vol(B_1tB_1)=\tr(\s_{\fT_y}|_{V^B}))/vol(BtB).\tag a$$
We can find a finite subset $L$ of $T(K)_0$ such that 
$BtB=\sqc_{\t\in L}B_1tB_1\t$. It follows that

(b) $vol(BtB)=vol(B_1tB_1)\sha(L)$
\nl
and $\s_{\fT_y}=\sum_{\t\in L}\s_{f_t}\s(\t)$ as linear maps $V@>>>V$. Restricting this equality
to $V^B$ and using the fact that $\s(\t)$ acts as identity on $V^B$ we obtain

(c) $\s_{\fT_y}|_{V^B}=\sha(L)\s_{f_t}|_{V^B}$ 
\nl
as linear maps $V^B@>>>V^B$. Clearly, (a) follows from (b) and (c). This completes the proof of Theorem 4.3.

\smallskip
The following result will not be used in the rest of the paper.
\proclaim{Proposition 5.5} If $y\in Y^{++}$ and $t\in T(K)_y$ then $BtB\sub{}^G\!T(K)_y$.
\endproclaim
\noindent
{\it Proof.} It is enough to show that $tz\sub{}^G\!T(K)_y$ for any $z\in B$. We can write $z=t_0z'$
where $t_0\in T(K)_0,z'\in B_1$. We have $tz=tt_0z'$ where $tt_0\in T(K)_y=T(K)^\spa_y$ (here we 
use that $y\in Y^{++}$). Using Lemma 5.3 we have $tt_0z'\in{}^G\!(tt_0T_1)\sub{}^G\!T(K)_y$.
This completes the proof. \qed

\subhead 5.6\endsubhead
In the remainder of this section we assume that $G$ is adjoint.
In this case the irreducible representations $(V,\s)$ as in 4.2 (up to isomorphism)
are known to be in bijection with the irreducible finite dimensional representations of the Hecke algebra $H$
(see \cite{BO}) by $(V,\s)\m V^B$. The irreducible finite dimensional representations of $H$ have been classified in \cite{KL1} in terms of geometric data.  
Moreover in \cite{L} an algorithm to compute the dimensions of the (generalized) weight spaces
of the action of the commutative semigroup $\{\fT_y;y\in Y^+\}$ on any tempered $H$ module is given. 
In particular the right hand side of the equality in Theorem 4.3 (hence also $\ph_V(t)$ in that Theorem)
is computable when $V$ is tempered.


\widestnumber\key{KL1}
\Refs
\ref\key{AK}\by J. Adler and J. Korman\paper The local character expansions near a tame, semisimple element\jour American J. of Math.\vol129\yr2007\pages381-403\endref
\ref\key{BO}\by A.Borel\paper Admissible representations of a semisimple group over a local field with fixed vectors under
an Iwahori subgroup\jour Inv.Math\yr1976\vol35\pages233-259\endref
\ref\key{C1}\by W.Casselman\paper The Steinberg character as a true character\jour Harmonic analysis on 
homogeneous spaces Proc.Symp.Pure Math.\vol26\yr1974\pages413-417\endref
\ref\key{C2}\by W.Casselman\paper Characters and Jacquet modules\jour Math.Ann.\vol230
\yr1977\pages101-105\endref
\ref\key{H}\by Harish-Chandra\paper Harmonic analysis on reductive
p-adic groups\inbook Harmonic analysis on homogeneous spaces\bookinfo Proc. Sympos.
Pure Math., Vol. XXVI \pages 167-192\publ Amer. Math. Soc., Providence, R.I.\yr1973\endref
\ref\key{IM}\by N.Iwahori and H.Matsumoto\paper On some Bruhat decomposition and the structure of the Hecke rings 
of $\fp$-adic Chevalley groups\jour Publ. Math\'ematiques IHES\vol25\yr1965\pages5-48\endref
\ref\key{KL1}\by D.Kazhdan and G.Lusztig\paper Proof of the Deligne-Langlands conjecture for Hecke algebras\jour
Inv.Math.\yr1987\vol87\pages153-215\endref
\ref\key{KL2}\bysame\paper Fixed point varieties on affine flag 
manifolds \jour Israel J.Math.\yr1988\vol62\pages129-168\endref
\ref\key{L}\by G.Lusztig\paper Graded Lie algebras and intersection cohomology
\jour Representation theory of algebraic groups and quantum groups, Progr. Math. Birkhauser, Springer \vol284\yr2010\pages191-224\endref
\ref\key{MA}\by H.Matsumoto\paper Fonctions sph\'eriques sur un groupe semi-simple $p$-adique\jour C.R.Acad.Sci. 
Paris\vol269\yr1969\pages829-832\endref
\ref\key{MP}\by A. Moy and G. Prasad \paper Unrefined minimal K-types for p-adic groups \jour Invent. Math. \vol 116 \yr 1994 no. 1-3\pages 393-408\endref
\ref\key{R}\by F.Rodier\paper Sur le caractere de Steinberg\jour Compositio Math.\vol59\yr1986\pages147-149\endref
\ref\key{S}\by J.A.Shalika\paper On the space of cusp forms on a $p$-adic Chevalley group\jour Ann.Math,\vol92
\yr1070\pages 262-278\endref
\ref\key{St}\by R. Steinberg\paper A geometric approach to the representations of the full linear group over a Galois field\jour Trans.Amer.Math.Soc.\vol71\yr1951\pages274–282\endref
\endRefs
\enddocument